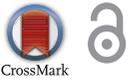

cogent
mathematics

# The plastic number and its generalized polynomial

Vasileios Iliopoulos[1]*

**Abstract:** The polynomial $X^3 - X - 1$ has a unique positive root known as plastic number, which is denoted by $\rho$ and is approximately equal to 1.32471795. In this note, we study the zeroes of the generalised polynomial $X^k - \sum_{j=0}^{k-2} X^j$, for $k \geq 3$, and prove that its unique positive root $\lambda_k$ tends to the golden ratio $\phi = \frac{1+\sqrt{5}}{2}$ as $k \to \infty$. We also derive bounds on $\lambda_k$ in terms of Fibonacci numbers.

**Subjects:** Science; Mathematics & Statistics; Advanced Mathematics; Number Theory; Pure Mathematics

## 1. Introduction

The recurrence $F_n = F_{n-1} + F_{n-2}$, with initial values $F_0 = 0$ and $F_1 = 1$ yields the celebrated Fibonacci numbers. It is well known that for $n \geq 0$

$$F_n = \frac{\phi^n - (1-\phi)^n}{\sqrt{5}},$$

where $\phi = \frac{1+\sqrt{5}}{2}$ is the positive root of the characteristic polynomial $X^2 - X - 1$, known as *golden ratio*.

One can readily generalise the recurrence and define the $k \geq 2$ order Fibonacci sequence $F_n = F_{n-1} + \cdots + F_{n-k}$, with initial conditions $F_0 = \cdots = F_{k-2} = 0$ and $F_{k-1} = 1$. The characteristic polynomial of this recurrence is $X^k - X^{k-1} - \cdots - X - 1$. Its zeroes are much studied in literature: we refer to Martin (2004), Miles (1960), Miller (1971), Wolfram (1998), and Zhu and Grossman (2009), where it is proved that the unique positive root tends to 2, as $k \to \infty$. Series representations for this root are derived in Hare, Prodinger, and Shallit (2014) by Lagrange inversion theorem.

### ABOUT THE AUTHOR

Vasileios Iliopoulos obtained his PhD from the Department of Mathematical Sciences of University of Essex in 2014 and is currently working as a part-time lecturer there. His research interests are on analytic and enumerative combinatorics, operational research and elementary number theory.

### PUBLIC INTEREST STATEMENT

In this paper, the author considers a result of Siegel (1944), who proved that the positive root of the polynomial $\mathcal{F}_k(X) = X^k - \frac{X^{k-1}-1}{X-1}$, where $k$ is an odd integer greater than or equal to 3, tends to the golden ratio $\frac{1+\sqrt{5}}{2}$, as $k \to \infty$.

The presented proof is elementary and simple. The author also obtains bounds on the positive zero of $\mathcal{F}_k(X)$ in terms of Fibonacci numbers.

cogent ·· mathematics

In this note, we turn our attention to the positive zero of the polynomial $X^3 - X - 1$, known as *plastic number*, which will throughout be denoted by $\rho$ and is equal to $\sqrt[3]{1 + \sqrt[3]{1 + \sqrt[3]{1 \ldots}}}$ (Finch, 2003). The plastic number was introduced by van der Laan (1960). The recurrence relation is $a_n = a_{n-2} + a_{n-3}$, with initial conditions $a_0 = a_1 = a_2 = 1$ and defines the integer sequence, known as Padovan sequence Stewart (1996). Although the bibliography regarding the analysis of Fibonacci numbers is quite extensive, it seems not to be this case regarding the plastic number.

In the next section, we examine a generalisation of the Padovan sequence and its associated characteristic polynomial and derive bounds on the unique positive root of the polynomial $X^k - X^{k-2} - \cdots - X - 1$. The presented study utilises results from the theory of linear recurrences and elementary Calculus, where the nature of roots of this polynomial is investigated. Furthermore, the bounds on the largest root are developed using identities of Fibonacci numbers.

## 2. The generalised sequence

Consider the recurrence

$$a_n = \sum_{l=2}^{k} a_{n-l}$$

for $k \geq 3$ and initial conditions $a_0 = \cdots = a_{k-1} = 1$. For $k = 3$, we obtain as a special case the Padovan sequence. A lemma follows regarding the roots of its characteristic polynomial.

LEMMA 2.1 *The polynomial $\mathcal{F}_k(X) = X^k - X^{k-2} - \cdots - X - 1$ has $k$ simple roots. If $k$ is odd, the polynomial has a unique real root $\lambda_k \in (1, \phi)$ and $k - 1$ complex roots. When $k$ is even, the roots of the polynomial are $\lambda_k \in (1, \phi)$, $-1$ and $k - 2$ complex zeroes.*

*Proof* It can be easily seen that neither $0$ nor $1$ are roots of $\mathcal{F}_k(X)$. Following Miles (1960) and Miller (1971), it is convenient to work with the polynomial

$$(X - 1)\mathcal{F}_k(X) = X^{k+1} - X^k - X^{k-1} + 1. \tag{1}$$

Differentiating Equation 1, we obtain

$$\left((X - 1)\mathcal{F}_k(X)\right)' = (k + 1)X^k - kX^{k-1} - (k - 1)X^{k-2}. \tag{2}$$

Equation 2 is $0$, at $X = 0$ or at the roots of the quadratic polynomial:

$$(k + 1)X^2 - kX - (k - 1). \tag{3}$$

Its discriminant can be easily computed to $\Delta = 5k^2 - 4 > 0$, for all $k \geq 3$ and the two real roots of polynomial of (3) are

$$\beta_{1,2}(k) = \frac{k \pm \sqrt{5k^2 - 4}}{2(k + 1)}. \tag{4}$$

We identify the real roots by elementary means. Note that $\mathcal{F}_k(1) = 2 - k < 0$ and $\mathcal{F}_k(\phi) = \phi$ and apply Descartes' rule of signs to Equation 1, there is a unique positive root $\lambda_k$ in $(1, \phi)$, and for $k$ even, the unique negative root of the polynomial is $-1$. Also, the polynomial is monic and by Gauss's lemma the root $\lambda_k$ is irrational.

Observe that $\lambda_k \neq \beta_2(k)$, since $\beta_2(k)$ is negative for all $k \geq 3$ and $\lambda_k \neq \beta_1(k)$. For if they were equal, then by Rolle's theorem there is at least one root $\alpha$ of Equation 3 in $(1, \beta_1(k))$, but $\beta_2(k) < 0 < 1$ and considering the fundamental theorem of Algebra, which states that every polynomial with complex coefficients and degree $k$ has $k$ complex roots with multiplicities, we arrive in contradiction. This shows that the polynomial $(X - 1)\mathcal{F}_k(X)$ has $(k + 1)$ simple roots. We complete the proof noting that $\mathcal{F}_k(X)$ and $(X - 1)\mathcal{F}_k(X)$ are positive and increasing for $X > \lambda_k$ and negative for $1 < X < \lambda_k$.

A direct consequence of Lemma 2.1 is

COROLLARY 2.2 *The polynomial $\mathcal{F}_k(X)$ is irreducible on the ring of integer numbers $\mathbb{Z}$ if and only if $k$ is odd.*

Further, it is easy to prove that all complex zeroes of the polynomial are inside the unit circle. The next Lemma is from Miles ([1960](#)) and Miller ([1971](#)).

LEMMA 2.3 *(Miles, [1960](#); Miller, [1971](#)) For all the complex zeroes $\mu$ of the polynomial $\mathcal{F}_k(X)$, it holds that $|\mu| < 1$.*

*Proof* Assume that there exists a complex $\mu$ (and hence $\overline{\mu}$), with $1 < |\mu| < \lambda_k$. We have that $(\mu - 1)\mathcal{F}_k(\mu) = 0$ and

$$|\mu^{k+1}| = |\mu^k + \mu^{k-1} - 1|. \tag{5}$$

Applying the triangle inequality to Equation 5, we deduce that

$$(|\mu| - 1)\mathcal{F}_k(|\mu|) < 0,$$

which contradicts Lemma 2.1. Assuming now that $|\mu| > \lambda_k$, we have

$$|\mu^k| = \left| \sum_{j=0}^{k-2} \mu^j \right| \leq \sum_{j=0}^{k-2} |\mu^j|,$$

which is equivalent to $\mathcal{F}_k(|\mu|) \leq 0$ and again we arrive in contradiction. Finally, by the same reasoning it can be easily proved that there is no complex zero $\mu$, with either $|\mu| = \lambda_k$ or $|\mu| = 1$.

Lemma 2.3 implies that the solution of the generalised recurrence can be approximated by

$$a_n \approx C\lambda_k^n, \tag{6}$$

with negligible error term. In Equation 6, $C$ is a constant to be determined by the solution of a linear system of the initial conditions.

We now consider, more carefully, Equation 4

$$\beta_{1,2}(k) = \frac{k \pm \sqrt{5k^2 - 4}}{2(k + 1)}.$$

Observe that $\beta_1(k) = \frac{k + \sqrt{5k^2 - 4}}{2(k+1)}$ is increasing and bounded sequence. Furthermore,

$$\lim_{k \to \infty} \frac{k + \sqrt{5k^2 - 4}}{2(k + 1)} = \frac{1}{2} + \sqrt{\frac{5}{4}} = \phi. \tag{7}$$

Also, $\beta_2(k) = \frac{k - \sqrt{5k^2 - 4}}{2(k+1)}$ is decreasing and bounded and

$$\lim_{k \to \infty} \frac{k - \sqrt{5k^2 - 4}}{2(k + 1)} = \frac{1}{2} - \sqrt{\frac{5}{4}} = 1 - \phi. \tag{8}$$

From Equations 7 and 8, we deduce that two of the critical points of Equation 1, (recall that these are $0$ with multiplicity $(k-2)$, $\beta_1(k)$ and $\beta_2(k)$), converge to $\phi$ and $1 - \phi$. A straightforward calculation can show that $\beta_1(k)$ are points of local minima of the function $(X-1)\mathcal{F}_k(X)$ to the interval $(1, \lambda_k)$, so $\beta_1(k) < \lambda_k < \phi$ for all $k \geq 3$ and by squeeze lemma we have that $\lim_{k \to \infty} \lambda_k = \phi$.

We remark that $\rho$ is a Pisot–Vijayaraghavan number, a real algebraic integer having modulus greater to $1$ where its conjugates lie inside the unit circle (Bertin, Decomps-Guilloux, Grandet-Hugot, Pathiaux-Delefosse, & Schreiber, 1992). These numbers are named after Pisot (1938) and Vijayaraghavan (1941), who independently studied them. Siegel (1944) considered several families of polynomials and showed that the plastic number is the smallest Pisot–Vijayaraghavan number. By Lemmas 2.1 and 2.3, the positive zeroes of the polynomial $X^k - \sum_{j=0}^{k-2} X^j$, where $k$ is odd, are Pisot–Vijayaraghavan numbers. In case that $k$ is even, the positive roots of the polynomial $X^k - \sum_{j=0}^{k-2} X^j$ are Salem numbers (Salem, 1945). This family of numbers is closely related to the set of Pisot–Vijayaraghavan numbers. They are positive algebraic integers with modulus greater than $1$, where its conjugates have modulus no greater than $1$ and at least one root has modulus equal to $1$.

We have proved that for all $k \geq 3$,

$$\frac{k + \sqrt{5k^2 - 4}}{2(k+1)} < \lambda_k < \phi. \tag{9}$$

Using the identity $5F_k^2 = L_k^2 - 4(-1)^k$ (Hoggatt, 1969, Section 5), we have that for $k = F_{2t+1}$

$$\lambda_{F_{2t+1}} > \frac{F_{2t+1} + L_{2t+1}}{2(F_{2t+1} + 1)} \tag{10}$$

where $L_n$ is the $n$th Lucas number, defined by $L_n = L_{n-1} + L_{n-2}$ for $n \geq 2$, with initial conditions $L_0 = 2$ and $L_1 = 1$. Lucas numbers obey the following closed form expression for $n \geq 0$, (Hoggatt, 1969)

$$L_n = \phi^n + (1 - \phi)^n.$$

Now inequality 10 becomes

$$\begin{aligned}\lambda_{F_{2t+1}} &> \frac{F_{2t+1} + L_{2t+1}}{2(F_{2t+1} + 1)} \\ &= \frac{F_{2(t+1)}}{F_{2t+1} + 1}.\end{aligned}$$

Actually, inequality 10 is valid when $\dfrac{k + \sqrt{5k^2 - 4}}{2(k+1)}$ is quadratic irrational. A stronger result is the following Theorem.

THEOREM 2.4   For $k \geq 3$, it holds that

$$\frac{F_{k+1}}{F_k + 1} < \lambda_k < \frac{F_{k+1}}{F_k}.$$

*Proof*   For $k = 3$, we have that

$$\frac{F_4}{F_3 + 1} < \rho < \frac{F_4}{F_3},$$

where $\lambda_3 := \rho$. Since for $k > 3$, $\frac{F_{k+1}}{F_k} > 1$ and $\frac{F_{k+1}}{F_k + 1} > 1$, by Lemma 2.1 it suffices to show that $\left(\frac{F_{k+1}}{F_k + 1} - 1\right)\mathcal{F}_k\left(\frac{F_{k+1}}{F_k + 1}\right) < 0$ and $\left(\frac{F_{k+1}}{F_k} - 1\right)\mathcal{F}_k\left(\frac{F_{k+1}}{F_k}\right) > 0$. Setting $X = \frac{F_{k+1}}{F_k}$ to Equation 1, we have to prove that

$$\left(\frac{F_{k+1}}{F_k}\right)^{k-1}\left(\left(\frac{F_{k+1}}{F_k}\right)^2 - \frac{F_{k+1}}{F_k} - 1\right) > -1.$$

The previous inequality is the same as

$$\left(\frac{F_{k+1}}{F_k}\right)^2 - \frac{F_{k+1}}{F_k} - 1 > -\left(\frac{F_k}{F_{k+1}}\right)^{k-1}. \tag{11}$$

The left-hand side of inequality 11 is

$$\frac{F_{k+1}^2 - F_k F_{k+1} - F_k^2}{F_k^2} = \frac{F_{k+1}(F_{k+1} - F_k) - F_k^2}{F_k^2}$$

$$= \frac{F_{k+1} F_{k-1} - F_k^2}{F_k^2}$$

$$= \frac{(-1)^k}{F_k^2},$$

by Cassini's identity (Hoggatt, 1969).

Thus, (11) is true for all $k$.

In order to prove that $\frac{F_{k+1}}{F_k+1} < \lambda_k$, we have to equivalently show that

$$\left(\frac{F_{k+1}}{F_k+1}\right)^2 - \frac{F_{k+1}}{F_k+1} - 1 < -\left(\frac{F_k+1}{F_{k+1}}\right)^{k-1}. \tag{12}$$

We then have

$$\left(\frac{F_{k+1}}{F_k+1}\right)^2 - \frac{F_{k+1}}{F_k+1} - 1 = \frac{F_{k+1}^2 - F_k F_{k+1} - F_{k+1} - F_k^2 - 2F_k - 1}{(F_k+1)^2}$$

$$= \frac{(-1)^k - (F_{k+1} + 2F_k + 1)}{(F_k+1)^2}$$

$$= \frac{(-1)^k - L_{k+1} - 1}{(F_k+1)^2},$$

which is

$$L_{k+1} F_{k+1}^{k-1} > (F_k+1)^{k+1}.$$

Using that

$$L_{k+1} = \phi^{k+1} + (1-\phi)^{k+1} \quad \text{and} \quad F_{k+1} = \frac{\phi^{k+1} - (1-\phi)^{k+1}}{\sqrt{5}}$$

the identity

$$F_{2(k+1)} > (F_k+1)^2$$

can be easily proven by induction, and thus completes the proof.

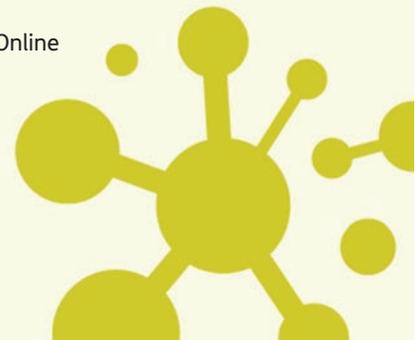

### Acknowledgements
I thank the anonymous referees for their helpful suggestions.

### Funding
The authors received no direct funding for this research.

### Author details
Vasileios Iliopoulos[1]

E-mail: viliop@essex.ac.uk

[1] Department of Mathematical Sciences, University of Essex, Colchester, UK.

### Citation information
Cite this article as: The plastic number and its generalized polynomial, Vasileios Iliopoulos, *Cogent Mathematics* (2015), 2: 1023123.

### References
Bertin, M. J., Decomps-Guilloux, A., Grandet-Hugot, M., Pathiaux-Delefosse, M. & Schreiber, J. P. (1992). *Pisot and Salem numbers*. Basel: Birkhäuser.

Finch, S. R. (2003). *Mathematical constants*. New York, NY: Cambridge University Press.

Hare, K., Prodinger, H., & Shallit, J. (2014). Three series for the generalized golden mean. *Fibonacci Quarterly, 52*, 307–313.

Hoggatt, V. E. Jr., (1969). *The Fibonacci and Lucas numbers*. Boston, MA: Houghton Mifflin.

Martin, P. A. (2004). The galois group of $x^n - x^{n-1} - \ldots - x - 1$. *Journal of Pure and Applied Algebra, 190*, 213–223.

Miles, E. P., Jr (1960). Generalized Fibonacci numbers and associated matrices. *American Mathematical Monthly, 67*, 745–752.

Miller, M. D. (1971). On generalized Fibonacci numbers. *American Mathematical Monthly, 78*, 1108–1109.

Pisot, C. (1938). La répartition modulo 1 et les nombres algébriques. *Annali di Pisa, 7*, 205–248.

Salem, R. (1945). Power series with integral coefficients. *Duke Mathematical Journal, 12*, 153–172.

Siegel, C. L. (1944). Algebraic numbers whose conjugates lie in the unit circle. *Duke Mathematical Journal, 11*, 597–602.

Stewart, I. (1996). Tales of a neglected number. *Scientific American, 274*, 102–103.

van der Laan, D. H. (1960). *Le Nombre Plastique: quinze leçons sur l'Ordonnance architectonique*. Leiden: Brill.

Vijayaraghavan, T. (1941). On the fractional parts of the powers of a number, II. *Proceedings of the Cambridge Philosophical Society, 37*, 349–357.

Wolfram, D. A. (1998). Solving generalized Fibonacci recurrences. *Fibonacci Quarterly, 36*, 129–145.

Zhu, X., & Grossman, G. (2009). Limits of zeros of polynomial sequences. *Journal of Computational Analysis and Applications., 11*, 140–158.